\documentclass[12pt]{article}
\usepackage{amsfonts}
\usepackage{amsmath}
\usepackage{pdfsync}

\newtheorem{theorem}{Theorem}
\newtheorem{definition}{Definition}
\newtheorem{lemma}{Lemma}

\newtheorem{corollary}{Corollary}
\begin{document}

\tableofcontents

%
%
%
%
%
%
%
%
%

\title[QMC. Unification approach]{Quantum Markov chains:\\  A unification approach}

\author[Luigi Accardi]{ Luigi Accardi }



\thanks{This work was completed with the support of our\TeX-pert.}
\author[Abdessatar Souissi]{Abdessatar Souissi}

\author[El Gheteb Soueidy]{El Gheteb Soueidy}




\textbf{abstract}
In the present paper we study a unified approach for Quan-
tum Markov Chains. A new quantum Markov property that generalizes
the old one, is discussed. We introduce Markov states and chains on
general local algebras, possessing a generic algebraic property, includ-
ing both Boson and Fermi algebras. The main result is a reconstruction
theorem for quantum Markov chains in the mentioned kind of local alge-
bras. Namely, this reconstruction allows the reproduction of all existing
examples of quantum Markov chains and states.

\section{Introduction and notations}\label{Introd}

Quantum Markov chains on infinite tensor product of matrix algebras were introduced
in \cite{[Ac74]} as a non--commutative analogue of classical Markov chains.
In \cite{[AcFr80]} the distinction between \textit{quantum Markov chains} and the subclass of
\textit{quantum Markov states} was introduced and a structure theorem for the latter class
was proved.
A sub--class of Markov chains, re--named finitely correlated states, was
shown to coincide with the so-called valence bond states introduced in the late 1980s
in the context of anti--ferromagnetic Heisenberg models (see \cite{F.N.W}).\\
In \cite{[AcFiMukh05]} the notion of quantum Markov chain was extended to states on the CAR algebra.
In \cite{AccMuSa} concrete models rising naturally from quantum statistical physics
were investigated in quantum spin algebras.\\

In the framework of more general $^*$--algebras a definition of Markov chains is still
missing. Namely, the following problems are still open
 \begin{itemize}
  \item An definition of Markov chains on $^*$--algebras more general than infinite
  tensor products of $^*$--algebras or CAR algebras.
 \item A reconstruction of Markov chains starting from the associated correlation functions.
 \end{itemize}
In this paper we solve the mentioned problems for an important class of quasi--local
$^*$--algebras for which the local algebras are linearly generated by "ordered products"
(see condition \eqref{A-gen-ord-prod} below).
These algebras include the infinite tensor products of type I factors and the Fermi algebra
generated by the canonical anti--commutation relations (CAR) (see \cite{BR} and \cite{BR2}).\\

\noindent The organization of the paper is the following.
In section 3 we introduce a formulation of the Markov property with respect to a backward
filtration $\{\mathcal A_{n]}\}$ that generalizes the Markov property introduced in \cite{[Acc75]}.
Section 4, is devoted to the definition of backward Markov states and chains in the
considered $*$--algebra $\mathcal A$ together with an existence theorem for Markov chains for given
sequence of boundary conditions.
In section 5, we state our main result which concerns a reconstruction of Markov chains starting
from a given sequence of transition expectations. We then prove that this result extends
the corresponding structure theorems in the tensor and the Fermi case.

\section{Notations and preliminaries}

Let $\mathcal{A}$ be a $*$--algebra and $\{\mathcal{A}_n\}_{n\in\mathbb N}$ a sequence
of its $*$--subalgebras. Unless otherwise specified, all $*$--algebras considered in the
following are complex unital, i.e. with identity.
For a given sub--set $I\subset \mathbb N$, denote
$$
\mathcal{A}_I: = \bigvee_{n\in I}\mathcal{A}_n
$$
the $*$--algebra generated by the family $(\mathcal{A}_n)_{n\in I}$. In these notations one has
$$
I\subseteq J \Rightarrow \mathcal{A}_I \subseteq \mathcal{A}_J
$$
If $I=[0,n]$, we denote  $\mathcal{A}_{n]}:=\mathcal{A}_{[0,n]}$.\\
If $I$ consists of a single element $n \in \mathbb N $ we write
$$
\mathcal{A}_n := \mathcal{A}_{\{n\}}
$$
The cone of positive elements of $\mathcal{A}_I$ will be denoted by $\mathcal{A}_{I}^+$.
We assume that the ordered products
\begin{equation}\label{A-gen-ord-prod}
a_0a_1\dots a_n \qquad;\qquad a_{j}\in\mathcal{A}_{j}, \ j\in\{1,\dots n\}, \ n\in\mathbb{N}
\end{equation}
linearly generate the algebra $\mathcal{A}$. This implies that any state $\varphi$ on
$\mathcal{A}$ is uniquely determined by its values on the products of the form
(\ref{A-gen-ord-prod}) and that
$$
\mathcal{A}=\mathcal{A}_{\mathbb N}
$$
For every integer $n\in\mathbb N^*$ denote by $M_n\equiv M(n,\mathbb  C)$ the
algebra of all complex $n\times n$ matrices.
Let $\mathcal{A}$ and $\mathcal{B}$ be two $*$--algebras.
\begin{definition}
Let $\mathcal{A}_{F}$ be a sub--$*$--algebra of $\mathcal{A}$.
A linear map $P$ from $\mathcal{A}_{F}$ into $\mathcal{B}$ is said to be
{\bf $n$--positive} ($n\in\mathbb N^*$) if $\forall\,b_1,\dots,b_n\in\mathcal{B}$, $\forall\,a_1,\dots,a_n\in\mathcal{A}_{F}$
\begin{equation}\label{QMCLV11)}
\sum^n_{j,k=1} b^*_j\,P(a^*_j\, a_k)b_k\geq 0
\end{equation}
$P$ is called {\bf completely positive\/} if \eqref{QMCLV11)} holds for all $n\in\mathbb N^*$.\\
If $\mathcal{C}\subseteq\mathcal{B}$ is a $*$--algebra and \eqref{QMCLV11)} holds for any
$n\in\mathbb N^*$ and any
$b_1,\dots,b_n\in \mathcal{C}$, $P$ is called $\mathcal{C}$--{\bf completely positive\/}.
\end{definition}
\begin{definition}
A linear map $E^0$ from $\mathcal{A}$ into $\mathcal{B}$ is called a
{\bf Umegaki conditional expectation} if:\\
\textbf{(CE1)} $E^0(a)\geq0$, if $a\geq0$; $a\in A$,\\
\textbf{(CE2)} $E^0(ba)=b$ $E(a)$; $b\in \hbox{Range}(E)$ \ , \ $a\in A$,\\
\textbf{(CE3)} $E^0(a^*)=E(a)^*$, $\forall a\in A$,\\
\textbf{(CE4)} $E^0(1)=1$,\\
\textbf{(CE5)} $E^0(a)\cdot E^0(a)^*\leq E^0(aa^*)$.\\
If such an $E^0$ exists, the algebra $\mathcal{B}$ is called \textit{expected}.
\end{definition}
\textbf{Remark}. If $E^0:\mathcal{A}\to \mathcal{B}\subseteq\mathcal{A}$ is a Umegaki conditional
expectation, (CE1) implies that $\|E(a)\|\leq\|a\|$ ($\forall a\in A$; (CE2) and (CE4) imply that
$E^0$ is a norm one projection onto its range which coincides with the set of its fixed points.
(CE2), (CE3) and (CE4) imply that $\hbox{Range}(E^0)$ is a $*$--algebra and that
$E^0:\mathcal{A}\to \hbox{Range}(E^0)$ is completely positive. In particular (CE5) follows from
(CE1)--(CE4).\\
\begin{definition}\label{df-QCE}
A \textbf{non--normalized quasi--conditional expectation} with respect to the
triplet of unital $\ast$-algebras $\mathcal{C} \subseteq \mathcal{B} \subseteq \mathcal{A}$
is a completely positive, linear $*$--map $E :\mathcal{A} \to \mathcal B$ such that
$E(1)\ne 1$ and
\begin{equation}\label{prop-QCEl}
E(ca)=cE(a)\quad,  \quad \forall a\in\mathcal{A} \ , \ \forall c\in\mathcal C
\end{equation}
If $E(1)=1$, $E$ is called a \textbf{quasi--conditional expectation}.\\
\end{definition}
\textbf{Remark}.
Any Umegaki conditional expectation $E$ from $\mathcal{A}$ into $\mathcal B$
satisfying \eqref{prop-QCEl} is a quasi--conditional expectation with respect to the
triplet $\mathcal{C} \subseteq \mathcal{B} \subseteq \mathcal{A}$.
\begin{lemma}\label{QCE-prop-of-CPmaps}
Let $P :\mathcal{A} \to \mathcal B$ be a completely positive map. Define these sets
\begin{equation}\label{prop-QCE2}
CE(P,l):=
\{c\in\mathcal{A} \ : \
P(ca)=cP(a) \hbox{ and }P(c^*a)=c^*P(a) \ ,  \  \forall a\in\mathcal{A} \}
\end{equation}
\begin{equation}\label{prop-QCE3}
CE(P,r):=
\{c\in\mathcal{A} \ : \
P(ac)=P(a)c \hbox{ and }P(ac^*)=P(a)c^* \ , \  \forall a\in\mathcal{A} \}
\end{equation}
Then both $CE(P,l)$ and $CE(P,r)$ are $*$--algebras and
\begin{equation}\label{CE(P)}
CE(P,l) = CE(P,r) =: CE(P)
\end{equation}
If $P$ is identity preserving,
$$
CE(P)\subseteq \hbox{Fix}(P):=\hbox{Fixed points of }P
$$
\end{lemma}
\textbf{Proof}.
It is clear that both $CE(P,l)$ and $CE(P,r)$ are algebra and they are closed under involution
by assumption. \eqref{CE(P)} follows from the identity
$$
P(ca)=cP(a)\iff
P(a^*c^*) = P((ca)^*) = (P(ca))^* = (cP(a))^* = (P(a))^* c^* = P(a^*) c^*
$$
Since $a\in\mathcal{A}$ and $c\in\mathcal CE(P,l)$ are arbitrary and $\hbox{Range}(P)$ is closed
under involution, this implies that the set \eqref{prop-QCE2} is equal to the set \eqref{prop-QCE3}.

\begin{lemma}\label{QCE-prop-of-CPmaps}
Let $E$ be a quasi--conditional expectation as in Definition \ref{df-QCE}. Then
there exists a $*$--sub--algebra $\mathcal C_{max}$ of $\hbox{Range}(E)$ maximal with respect
to property \eqref{prop-QCEl} and such that $\mathcal C_{max}$.
\begin{equation}\label{Cmax=Fix(E)}
\mathcal C\subseteq\mathcal C_{max} \subseteq \hbox{Fix}(E)\subseteq\mathcal B
\end{equation}
\end{lemma}
\textbf{Proof}.
From Zorn Lemma it follows that there exists a $*$--sub--algebra $\mathcal C_{max}$ of
$\hbox{Range}(E)$ maximal with respect to property \eqref{prop-QCEl} and such that
$\mathcal C\subseteq\mathcal C_{max}$.
\eqref{Cmax=Fix(E)} then follows because we have seen that property \eqref{prop-QCEl} implies that
$\mathcal C_{max}\subseteq\hbox{Fix}(E)$.\\

\noindent {\bf Remark}.
Suppose that the algebra $\mathcal C_{max}$ in Lemma \ref{QCE-prop-of-CPmaps} is expected
and let $E^0:\mathcal A\to \mathcal C_{max}$ be a surjective Umegaki conditional expectation.
Any $K\in\mathcal A$ such that
$$
E^0(K^*K) = 1
$$
is called an $E^0$--\textbf{conditional amplitude}.
Denoting, for any sub--$*$--algebra $\mathcal{B}\subseteq\mathcal{A}$
$$
\mathcal{B}^{'}:=
\{a\in\mathcal{A} \ : \ a b= ba \ , \ \forall b\in\mathcal{B} \}
$$
the commutant of $\mathcal{B}$ in $\mathcal{A}$,  If $K\in\mathcal C'$ then the map
$$
E^0(K^*( \ \cdot \ ) K) : \mathcal{A} \to \mathcal{B}
$$
is a quasi--conditional expectation with respect to the
triplet $\mathcal{C} \subseteq \mathcal{B} \subseteq \mathcal{A}$.\\

\noindent {\bf Remark}. Every quasi--conditional expectation with respect to the triplet
$\mathcal{C} \subseteq \mathcal{B} \subseteq \mathcal{A}$ satisfies the conditions
\begin{equation}\label{prop-QCEr}
E(ac)=E(a)c \quad; \quad  a\in\mathcal A \ , \  c\in \mathcal C
\end{equation}
\begin{equation}\label{old_MP}
E(\mathcal C'\cap \mathcal{A})\subseteq \mathcal C'\cap B
\end{equation}

\section{A new formulation of the backward quantum Markov property}

\begin{definition}
A map $E$ from  $\mathcal{A}_{n+1]}$ into $\mathcal{A}_{n]} $ is said to enjoy the
\textbf{Markov property} with respect to the triplet
$\mathcal{A}_{n-1]}\subseteq \mathcal{A}_{n]}\subseteq \mathcal{A}_{n+1]} $ if
\begin{equation}\label{MP}
E(\mathcal{A}_{[n,n+1]}) \subseteq \mathcal{A}_{n-1]}^{'}\cap \mathcal{A}_{n}
\end{equation}
\end{definition}
{\bf Remark}.
In \cite{[Ac74]} it was claimed that the relation (\ref{old_MP}) can be considered as
a non--commutative formulation of the Markov property and it was shown that this claim is plausible
in the tensor case in which $A_{n]}=\bigotimes_{[0,n]} M_d(\mathbb C)$ for each $n$.
In this case in fact on has
\begin{equation}\label{old_MP2}
\mathcal{A}^{'}_{n-1]}\cap \mathcal{A}_{n+1]}=\mathcal{A}_{[n,n+1]}
\end{equation}
However in the Fermi case \eqref{old_MP2} is not satisfied, while our Definition
(\ref{MP}) applies to both cases (see section \ref{Fermi-case}).\\
{\bf Remark}.
From \eqref{prop-QCEl}, \eqref{prop-QCEr} and \eqref{MP}, it follows that for any $a_{n-1]}\in\mathcal{A}_{n-1]}$
and\\ $a_{[n,n+1]}\in\mathcal{A}_{[n,n+1]}$ one has
$$
E(a_{n-1]}a_{[n,n+1]})
=a_{n-1]}E(a_{[n,n+1]})
=E(a_{[n,n+1]})a_{n-1]}
=E(a_{[n,n+1]}a_{n-1]})
$$
Therefore \textbf{any Markov quasi--conditional expectation} $E_{n]}$ with respect
to the triplet $ \mathcal{A}_{n-1]}\subset \mathcal A_{n]}\subset \mathcal A_{n+1]}$ must
satisfy the following trace--like property
\begin{equation}\label{star_condition}
E_{n]}(ab)=E_{n]}(ba)  \quad ;\quad a \in \mathcal A_{n-1]},\  b \in \mathcal{A}_{[n,n+1]}
\end{equation}
\begin{definition}\label{Fwrd_TE}
A \textbf{backward Markov transition expectation} from\\
$\mathcal{A}_{n}\vee  \mathcal{A}_{n+1}$ to $\mathcal{A}_{n}$
is a completely positive identity preserving map
$$
E_{[n,n+1]}: \ \mathcal{A}_{n}\vee  \mathcal{A}_{n+1} \rightarrow \mathcal{A}_{n}
$$
satisfying the \textbf{Markov property} \eqref{MP}.\\

If $E_{[n+1,n]}$ is not identity preserving, we say that it is a non--normalized backward Markov
transition expectation.
\end{definition}
{\bf Remark}.
Any Markov quasi--conditional expectation $E_{n]}$ with respect to the triplet
$ \mathcal{A}_{n-1]}\subset \mathcal A_{n]}\subset \mathcal A_{n+1]}$
defines, by restriction to $\mathcal{A}_{[n,n+1]} $ a backward Markov transition expectation
$E_{[n+1,n]}$ from $\mathcal{A}_{n}\vee  \mathcal{A}_{n+1}$ to $\mathcal{A}_{n}$ with respect
to the same triplet.\\
We will prove in the section (\ref{recons-sec}) that any backward Markov transition expectation
$E_{[n+1,n]}$ from $\mathcal{A}_{n}\vee  \mathcal{A}_{n+1}$ to $\mathcal{A}_{n}$ with respect
to the triplet $ \mathcal{A}_{n-1]}\subset \mathcal A_{n]}\subset \mathcal A_{n+1]}$
arises in this way. To this goal we recall some properties of the non--commutative
Schur multiplication.
\begin{definition}  Let $\mathcal M$ be a $\ast$--algebra and let
$A=(a_{ij}), B=(b_{ij})\in M_n(\mathcal M)$.\\
The \textbf{Schur} product $A\circ B$ is defined by
\begin{equation}\label{tensor-schur}
A\circ B := (a_{ij}b_{ij})\in  M_n(\mathcal M)
\end{equation}
\end{definition}
\textbf{Remark}.   Note that  $A\circ B = B\circ A$ if and only if for each $i,j$  the elements $a_{ij}$ and $b_{ij}$ commute.
\begin{lemma}\label{schur}
Let $\mathcal M$ be a $*$-algebra and $\mathcal A$, $\mathcal B$ commuting sub--$*$--algebras
of $\mathcal M$. Then the Schur multiplication
$M_n(\mathcal A) \times M_n(\mathcal B) \mapsto M_n( \mathcal A \vee \mathcal B )$ is a positive map.
\end{lemma}
\textbf{Proof}.
Recall that by definition $A\in M_n(\mathcal A)$ is positive if and only if
it is a sum of elements of the form $A = C^\ast C$ with $C\in M_n(\mathcal A)$.
By linearity it will be sufficient to consider only elements of the form $A = C^\ast C$, i.e.
$$
a_{ij} = \sum_{h}c^\ast_{hi}c_{hj}\quad i,j\in \{1,n\}
$$
Let $A=C^\ast C\in M_n(\mathcal A)^{+}$, $B=D^\ast D\in M_n(\mathcal B)^{+}$.\\
For $X= (x_1, \cdots, x_n)^T \in (\mathcal A \vee \mathcal B )^{n}$, one has
$$
X^* A\circ B X=\sum_{i,j}x_{i}^\ast a_{ij}b_{ij}x_j\
=\sum_{i,j}x_{i}^\ast\left(\sum_{h}c^\ast_{hi}c_{hj}\right)\left(\sum_{k}d^\ast_{ki}d_{kj}\right)x_j
$$
$$
=\sum_{h,k}\sum_{i,j}x_{i}^\ast c^\ast_{hi}c_{hj}d^\ast_{ki}d_{kj}x_j
=\sum_{h,k}\sum_{i,j} (d_{ki}c_{hi}x_i)^\ast(d_{kj}c_{hj}x_j)
=\sum_{h,k}\bigl(\sum_{i} d_{ki}c_{ki}x_i\bigr)^\ast\bigl(\sum_{i}d_{ki}c_{hi}x_i\bigr)
$$
$=\sum_{h,k}\bigl|\sum_{i} d_{ki}c_{ki}x_i\bigr|^2\geq 0$.
Then $A \circ B\in M_n(\mathcal A \vee \mathcal B)_+$.\\

\begin{definition}\label{df-tensor-schur}
Let be given  two unital $\ast$--algebras $\mathcal M$ and $\mathcal V$.
If $A=[a_{ij}]\in M_n(\mathcal M)$ and $B=[b_{ij}]\in M_n(\mathcal V)$ their
\textbf{Schur tensor product} is defined by
\begin{equation}\label{tensor-schur}
A\circ^\otimes B := [a_{ij}\otimes b_{ij}]\in  M_n(\mathcal M \otimes \mathcal V)
\end{equation}
where $\otimes$ is the algebraic tensor product.
\end{definition}
\begin{lemma}\label{Tensor_Schur_lemma}
In the notations of Definition \ref{df-tensor-schur}, if $A$ and $B$ are positive then
$A\circ^\otimes B$ is also positive.
\end{lemma}
\textbf{Proof}. See \cite{Sumesh2016}.\\

\noindent We will use the following corollary of Lemma \ref{Tensor_Schur_lemma}.
\begin{corollary}\label{Cor-Tens-Schur-lm}
In the notations of Definition \ref{df-tensor-schur}, let $\mathcal C$ and $\mathcal D$
be mutually commuting sub--algebras of a $*$--algebra $\mathcal A$ and let
$u:\mathcal M\to \mathcal C$ be a $*$--homomorphism and $P:\mathcal N\to \mathcal D$
a completely positive map.
Then the map
$$
u\otimes P : m\otimes n \in\mathcal M\otimes\mathcal N \to u(m)P(n)
\in\mathcal C \vee \mathcal D
$$
is $\mathcal{A}_{n-1]}$--completely positive.
\end{corollary}
\textbf{Proof}. We have to prove that for each $n\in \mathbb{N}$ the map
$$
\sum_{i,k=1}^{n} m_{j}^*m_{k}\otimes n_{j}^*n_{k} \mapsto
\sum_{i,k=1}^{n}  u(m_{i}^*m_{k})P(n_{i}^*n_{k})
$$
is positive. Since $(P(n_{i}^*n_{k}))$ is positive because $P$ is completely positive and
$ u(m_{i}^*m_{k})$ is positive because $u$ is a $*$--homomorphism,
the thesis follows from Lemma \eqref{schur}.

\section{Backward Markov states and chains}\label{Backw-MCC}

\subsection{Backward Markov states }\label{Backw-MS}

\begin{definition}{\rm
A state $\varphi$ on $\mathcal{A}$ is said to be a \textbf{backward quantum Markov state}
if for every $n\in\mathbb{N}$ there exists a, non necessarily normalized, Markov quasi--conditional
expectation $E_{n]}$ with respect to the triplet
$\mathcal{A}_{n-1]}\subseteq\mathcal{A}_{n]}\subseteq\mathcal{A}_{n+1]}$ satisfying
\begin{equation}\label{QMS_def}
\varphi=\varphi \circ E_{n]}
\end{equation}
for each ordered product $a_0a_1\cdots a_n$ with $a_k\in\mathcal{A}_k, k=1,\cdots,n$.
}\end{definition}
\begin{theorem}\label{struc-MS}
Any Markov state $\varphi$ on $\mathcal{A}$ defines a pair $\{\varphi_0 \ , \ (E_{[n,n+1]})\}$
such that:\\
(i) For every $n\in\mathbb N$,
\begin{equation}\label{state-cond}
\varphi_0(E_{0]}(E_{1]}((\cdots E_{n-1]}(E_{n]}(1_{n+1}))\cdots))))  = 1
\end{equation}
(ii) $\forall\,n\in\mathbb N$,
$E_{[n,n+1]}:\mathcal{A}_{[n,n+1]}\to \mathcal{A}_{n-1]}^{'}\cap \mathcal{A}_{n} $ is a linear
completely positive map;\\
(iii) For every $n\in\mathbb N$, $a_i \in \mathcal{A}_i$, $0\leq i\leq n$,
\begin{equation}\label{5.2.1}
\varphi(a_0a_1\cdots a_n)
=\varphi_0(E_{0]}(a_0E_{1]}(a_1(\cdots E_{n-1]}(a_{n-1}E_{n]}(a_n))\cdots))))
\end{equation}
Conversely, given a pair $\{\varphi_0 \ , \ (E_{[n,n+1]})\}$ satisfying (i) and (ii) above,
for every $n\in\mathbb{N}$ there is a unique state $\varphi_{[0,n]}$ on $\mathcal{A}_{[0,n]}$
satisfying
\begin{equation}\label{5.2.1b}
\varphi_{[0,n]}(a_0a_1\cdots a_n)
=\varphi_0(E_{[0,1]}(a_0E_{[1,2]}(a_1(\cdots E_{[n-1,n]}(a_{n-1}E_{[n,n+1]}(a_n))\cdots)
\end{equation}
If the family of states $(\varphi_{[0,n]})$ is projective, in the sense that
\begin{equation}\label{proj-cond-QMC}
\varphi_{[0,n+1]}\Big|_{\mathcal{A}_{[0,n]}} = \varphi_{[0,n]} \qquad;\qquad\forall n\in\mathbb{N}
\end{equation}
then it defines a unique state $\varphi$ on $\mathcal{A}$.\\ $\varphi$ is a Markov state
if and only if the compatibility condition
\begin{equation}\label{comp-cond-MS}
\varphi_{[0,n]}\left(a_{n-1]}E_{[n,n+1]}(a_{n-1}a_{n+1})\right)
=\varphi_{[0,n]}\left(a_{n-1]}E_{[n,n+1]}(a_{n-1}E_{[n+2,n+1]}(a_{n+1}))\right)
\end{equation}
is satisfied for any $a_{n-1]}\in\mathcal{A}_{[0,n-1]}$, $a_{n}\in\mathcal{A}_{,n}$ and
$a_{n+1}\in\mathcal{A}_{n+1}$.
\end{theorem}
\textbf{Proof}. \textbf{Necessity}.
Let $\varphi$ be a Markov state on $\mathcal{A}$ and let $(E_{n]})$ denote the associated
sequence of Markov quasi--conditional expectations. The map
\begin{equation}\label{E[n,n+1]restr}
E_{[n,n+1]}:=E_{n+1]}\Big|_{\mathcal{A}_{[n,n+1]}}:=\hbox{restriction of $E_{n]}$ on }
\mathcal{A}_{[n,n+1]}
\end{equation}
satisfies condition (ii) being the restriction of a map satisfying it. Denote
$$
\varphi_{0}:=\varphi\Big|_{\mathcal{A}_{0}}:=\hbox{restriction of $\varphi$ on }
\mathcal{A}_{0}
$$
Then iterated application of (\ref{QMS_def}) leads to
$$
\varphi(a_0a_1\dots a_n)=\varphi(a_0\dots a_{n-1}E_{n]}(a_n))
=\varphi(a_0\dots a_{n-2}E_{n-1]}(a_{n-1}E_{n]}(a_n))\\
=\cdots
$$
$$
=\varphi_0(E_{0]}(a_0E_{1]}(a_1(\cdots E_{n-1]}(a_{n-1}E_{n]}(a_n))\cdots)
$$
which, due to \eqref{E[n,n+1]restr}, is equivalent to \eqref{5.2.1}. Finally
condition (i) is satisfied because $\varphi$ is a state.\\
\textbf{Sufficiency}. Let $\{\varphi_0 \ , \ (E_{[n,n+1]})\}$ be a pair satisfying (i) and (ii) above
and, for each $n\in\mathbb{N}$, let $E_{n]}$ be the unique Markov quasi--conditional expectation
with respect to the triplet $ \mathcal{A}_{n-1]}\subset \mathcal A_{n]}\subset \mathcal A_{n+1]}$
associated to $E_{[n+1,n]}$ according to Theorem \ref{general-QCE}.
Then the composition
$$
E_{0]}\cdots E_{n]}E_{n+1]}
$$
is a completely positive map.
From positivity and condition \eqref{state-cond} it follows that the linear functional
$$
\varphi_{[0,n]} :=
\varphi_{0}E_{0]}\cdots E_{n-1]}E_{n]}
$$
is a state on $\mathcal{A}_{n+1]}$ which by construction satisfies \eqref{5.2.1b}.\\
It is known that the projectivity of the family of states $(\varphi_{[0,n]})$ is equivalent to
the existence of a unique state $\varphi$ on $\mathcal{A}$ whose restriction on each
$\mathcal{A}_{n]}$ is equal to $\varphi_{[0,n]}$.
This state will be A Markov state if and only if condition \eqref{QMS_def}
is satisfied and this is equivalent to
$$
\varphi_{[0,n+1]}(a_{n-1]} a_n a_{n+1}) = \varphi\circ E_{n]} (a_{n-1]} a_na_{n+1})
=\varphi_{[0,n]} (a_{n-1]} E_{n]}(a_{n-1}\cdot a_{n+1}))
$$
$$
=\varphi_{[0,n+1]}(a_{n-1]} a_na_{n+1}\cdot 1_{n+2})
=\varphi (a_{n-1]} a_{n-1}E_{n+1]}(a_n\cdot 1_{n+2}))
$$
$$
=\varphi\circ E_{[0,n]} (a_{n-1]} a_{n}E_{n+1]}(a_n\cdot 1_{n+2}))
= \varphi_{[0,n]}\left(a_{n-1]}E_{[n,n+1]}(a_{n}E_{[n+2,n+1]}(a_{n+1}))\right)
=
$$
is satisfied for any $a_{n-1]}\in\mathcal{A}_{[0,n-1]}$, $a_{n}\in\mathcal{A}_{,n}$ and
$a_{n+1}\in\mathcal{A}_{n+1}$, which is \eqref{comp-cond-MS}.

\subsection {Backward Markov chains }\label{NN-Backw-MC}

We have seen that any Markov state $\varphi$ on $\mathcal{A}$ defines a pair
$\{\varphi_0 \ , \ (E_{[n,n+1]})\}$ satisfying conditions (i) and (ii) of Theorem \ref{struc-MS}.
However not every pair satisfying these two conditions defines a Markov state on $\mathcal{A}$:
this is the case if and only if the compatibility condition \eqref{comp-cond-MS} is satisfied.
However it can happen that the pair $\{\varphi_0 \ , \ (E_{[n,n+1]})\}$ defines through
formula \eqref{5.2.1b} a family of states $(\varphi_{[0,n]})$ with the property that the limit
\begin{equation}\label{lim-varphi[0,n]}
\lim_{N\to\infty} \varphi_{[0,n]} =: \varphi
\end{equation}
exists point--wise on $\mathcal{A}$. Since we know from Theorem \ref{struc-MS} that each
$\varphi_{[0,n]}$ is a state on $\mathcal{A}$, the same will be true for $\varphi$.
The class of states defined by \eqref{lim-varphi[0,n]} turned out to have several interesting
applications in the theory of quantum spin system (where only algebras of the form
$\bigotimes_{n\in V}$ are considered, $V$ being the set of vertices of a Cayley tree).
If the limit \eqref{lim-varphi[0,n]} exists, because of assumption \eqref{A-gen-ord-prod}
it is uniquely determined by its values on the products of the form \eqref{A-gen-ord-prod}.
Therefore, because of \eqref{5.2.1b},  the limit \eqref{lim-varphi[0,n]} exists if and only if the
limit
\begin{equation}\label{lim-varphi[0,n]b}
\lim_{k\to\infty}\varphi_{[0,n+k]}(a_0a_1\cdots a_n1_{n+1}\cdots 1_{n+k})
=\lim_{k\to\infty}
\end{equation}
$$
\varphi_0(E_{[0,1]}(a_0E_{[1,2]}(a_1(\cdots E_{[n,n+1]}(a_{n}
E_{[n+1,n+2]}\cdots E_{[n+k-1,n+k]}(1_{n+k}))\cdots)
$$
exists for any $n\in\mathbb{N}$ and any $a_{j}\in\mathcal{A}_{j}$, $j\in\{1,\dots,n\}$. \\
Notice that, if the pair $\{\varphi_0 \ , \ (E_{[n,n+1]})\}$ satisfies conditions (i) and (ii) of Theorem \ref{struc-MS}, then
\begin{equation}\label{df-bn}
b_{n} := E_{[n,n+1]}(1_{n+1})\in\left(\mathcal{A}_{n-1]}^{'}\cap\mathcal{A}_{n}\right)_+
\qquad;\qquad \forall n\in\mathbb{N}
\end{equation}
It is clear that, if the sequence $(b_{n})$ defined by \eqref{df-bn} satisfies these condition
\begin{equation}\label{martingale_equation1}
E_{[n,n+1]}(b_{n+1})=b_{n}
\end{equation}
(see \cite{[Ac74]}, Lemma 1 for the tensor analogue of this condition) then for any $k\ge 2$
$$
\varphi_0(E_{[0,1]}(a_0E_{[1,2]}(a_1(\cdots E_{[n,n+1]}(a_{n}
E_{[n+1,n+2]}\cdots E_{[n+k-2,n+k-1]}(b_{n+k-1})\cdots)))))
$$
$$
\varphi_0(E_{[0,1]}(a_0E_{[1,2]}(a_1(\cdots E_{[n,n+1]}(a_{n}
E_{[n+1,n+2]}\cdots E_{[n+k-3,n+k-2]}(b_{n+k-2})\cdots)))))
$$
$$
=\dots =
\varphi_0(E_{[0,1]}(a_0E_{[1,2]}(a_1(\cdots E_{[n,n+1]}(a_{n}b_{n+1})))))
$$
i.e. the sequence $(\varphi_0(E_{[0,1]}(a_0E_{[1,2]}(a_1(\cdots E_{[n,n+1]}(a_{n}
E_{[n+1,n+2]}\cdots E_{[n+k-1,n+k]}(1_{n+k}))\cdots))_{k\ge 2}$  is constant, hence the
limit \eqref{lim-varphi[0,n]b} exists trivially and is equal to
\begin{equation}\label{lim-varphi[0,n]c}
\lim_{k\to\infty}\varphi_{[0,n+k]}(a_0a_1\cdots a_n1_{n+1}\cdots 1_{n+k})=
\end{equation}
$$
=\varphi_0(E_{[0,1]}(a_0E_{[1,2]}(a_1(\cdots E_{[n,n+1]}(a_{n}b_{n+1})\cdots))))
$$
\textbf{Remark}.
Equation (\ref{martingale_equation1}) means that the sequence $(b_{n})$ is a $(E_{n]})$--martingale.\\
\textbf{Remark}. Condition \eqref{martingale_equation1} is only sufficient to guarantee the
existence of the limit \eqref{lim-varphi[0,n]b}. Moreover, if  $\mathcal{A}$ is a $C^*$--algebra,
using the compactness of the states on $\mathcal{A}$ one can show that there is always at least
one sub--sequence of $(\varphi_{[0,n]})$ which defines a state on $\mathcal{A}$.
This justifies the following definition.
\begin{definition}\label{def_Backw_MC}
Let $\{\varphi_0 \ , \ (E_{[n,n+1]})\}$ be a pair satisfying conditions (i) and (ii) of
Theorem \ref{struc-MS} and let $(b_n)_{n\geq 0}\}$ be a sequnce of positive elements
$b_n\in \mathcal{A}_n$. Any state $\varphi$ on $\mathcal{A}$ satisfying
\begin{equation}\label{lim-varphi[0,n]d}
\varphi(a_0a_1\cdots a_n)=
\end{equation}
$$
=\lim_{k\to\infty}\varphi_0(E_{[0,1]}(a_0E_{[1,2]}(a_1(\cdots E_{[n,n+1]}(a_{n}
E_{[n+1,n+2]}\cdots E_{[n+k-1,n+k]}(b_{n+k+1}))\cdots)
$$
for any $n\in\mathbb{N}$ and any $a_{j}\in\mathcal{A}_{j}$, $j\in\{1,\dots,n\}$ is
called a {\bf  backward Markov chain} on $\mathcal{A}$ and the sequence $(b_{n})_{n\ge 0}$
is called the sequence of \textbf{boundary conditions} with respect to $(E_{n]})_{n\ge 0}$.
\end{definition}

\begin{theorem}\label{Suff_Cond-thm}{\rm
A sufficient condition for a triplet $\{\varphi_0, (E_{n]}), (b_n)\}$ to define a backward Markov chain is the existence of $c_{n} \in \mathcal{A}_{n-1]}^{'}$ for each $n$ such that
\begin{equation}\label{commutation_condition}
b_{n}= c_{n}^*c_{n}
\end{equation}
\begin{equation}\label{initial_cond}
\varphi_0(b_0) = 1
\end{equation}
\begin{equation}\label{martingale_equation}
E_{n]}(b_{n+1})=b_{n}
\end{equation}
Moreover, under these conditions the limit  exists in the strongly finite sense.\\
}\end{theorem}
\textbf{Proof}. {\rm
Using (\ref{martingale_equation}) one gets for every ordered product $a_0 a_1 \cdots a_n \in \mathcal{A}_{n]}$
\begin{eqnarray*}
&&\varphi_{0}(E_{0]}(a_0{E}_{1]}
(a_1(\cdots{E}_{n-1]}(a_{n-1}{E}_{n]}(a_n({E}_{n+1]}(1_{n+1}{E}_{n+2]}(\cdots {E}_{n+k]}(b_{n+k+1})) \cdots))))\cdots)))) \\
&&=\varphi_{0}(E_{0]}(a_0{E}_{1]}(a_1(\cdots
{E}_{n-1]}(a_{n-1}{E}_{n]}(a_n({E}_{n+1]}(1_{n+1}{E}_{n+2]}(\cdots {E}_{n+k-1]}(1b_{n+k})  \cdots)))) \cdots))))\\
&&\vdots\\
&&=\varphi_{0}(E_{0]}(a_0{E}_{1]}(a_1(\cdots {E}_{n-1]}(a_{n-1}{E}_{n]}(a_n {E}_{n+1]}( b_{n+2})))\cdots))))\\
&&=\varphi_{0}(E_{0]}(a_0{E}_{1]}(a_1(\cdots {E}_{n-1]}(a_{n-1}{E}_{n]}(a_n b_{n+1}))\cdots))))\\
\end{eqnarray*}
Then, the limit in the right hand side of (\ref{lim-varphi[0,n]d}) stabilizes at $n+1$,
i.e. it is equal to
$$
\varphi(a_0a_1\cdots a_n)
=\varphi_{0}(E_{0]}(a_0{E}_{1]}(a_1(\cdots {E}_{n-1]}(a_{n-1}{E}_{n]}(a_n b_{n+1}))\cdots))))
$$
\begin{equation}\label{finite_vol_BwrdMC}
=\varphi_{0}\circ E_{0]}\circ E_{1]}\cdots \circ E_{n]}(a_0a_1\cdots a_n b_{n+1})
\end{equation}
Now from (\ref{commutation_condition}) one gets
$$
E_{n]}(a_{n]}b_{n+1}) = E_{n]}(c_{n+1}^* a_{n]}c_{n+1})
$$
Therefore, the map
$$
E_{n],b}: a\in \mathcal{A}_{n]} \mapsto E_{n]}(a b_{n+1})\in  \mathcal{A}_{n]}
$$
is completely positive as a composition of completely positive maps.\\
Then through (\ref{finite_vol_BwrdMC}), the functional $\varphi$ is positive.\\
Therefore, taking into account (\ref{initial_cond}) one obtain $\varphi$ is
a quantum Markov chain in the sense of Definition \ref{def_Backw_MC}.

\section{Reconstruction theorem for backward Markov chain}\label{recons-sec}

Since the $*$-algebra $\mathcal A$ is linearly generated by ordered products of the form (\ref{A-gen-ord-prod}). Then via Zorn's lemma it admits a linear  basis which consists only of such ordered products.\\
We deal with the case where $\mathcal A$ has the following property: \\
for every
$m\in\mathbb N$, $B_m = \{e^{(m)}_{i_m}\}_{i_m\in I_m}$ is a linear basis of
$\mathcal A_m$, then
\begin{equation}\label{basis2}
B_{n]} := \{e_{i_0}^{(0)}e_{i_1}^{(1)}\cdots e_{i_n}^{(n)}\quad ; \quad (i_0,\cdots ,i_n)\in I_0\times \cdots \times  I_n\}
\end{equation}
is a linear basis of the $*$-algebra $\mathcal A_{n]}$, for each $n$ and
$$
B_{n+1]} = \{ e_{i_{n]}}^{n]}e_{i_{n+1}}^{(n+1)}\quad;\quad i_{n]}\in I_{n]}:= I_0\times \cdots \times  I_n  \ \, i_{n+1}\in I_{n+1}\}
$$

Let be given a Umegaki conditional expectation $E_{[n+1,n]}^{0}$ from $\mathcal{A}_{[n,n+1]}$ into
$\mathcal{A}_{n-1]}^{'} \cap \mathcal{A}_{[n,n+1]}$
and due to \eqref{basis2}, we can define
$$
\widetilde E_{n]}^{0}: \ \mathcal{A}_{n-1]} \vee \mathcal{A}_{[n,n+1]} \mapsto \mathcal{A}_{n+1]}
$$
as linear extension of
$$
\widetilde E_{n]}^{0}\left(e_{i_0}^{(0)}e_{i_1}^{(1)}\cdots e_{i_{n-1}}^{(n-1)}e_{i_{n}}^{(n)}e_{i_{n+1}}^{(n+1)}\right):
= e_{i_0}^{(0)}e_{i_1}^{(1)}\cdots e_{i_{n-1}}^{(n-1)}E_{[n+1,n]}^{0}\left(e_{i_{n}}^{(n)}e_{i_{n+1}}^{(n+1)}\right)
$$
which satisfies

\begin{equation}\label{Umegaki-cond}
\widetilde E_{n]}^{0}\left(e_{i_{n}}^{(n)}e_{i_{n+1}}^{(n+1)}e_{i_0}^{(0)}e_{i_1}^{(1)}\cdots e_{i_{n-1}}^{(n-1)}\right)
\end{equation}
$$
=E_{[n+1,n]}^{0}\left(e_{i_{n}}^{(n)}e_{i_{n+1}}^{(n+1)}\right)(e_{i_0}^{(0)}e_{i_1}^{(1)}\cdots e_{i_{n-1}}^{(n-1)})
$$
\textbf{Remark}.\rm{ One can see that from \eqref{Umegaki-cond}, we obtain
$$
\widetilde E_{n]}^{0}(ab)=\widetilde E_{n]}^{0}(ba),\ \ \text{for each} \  a \in \mathcal{A}_{n-1]},\ \  b\in \mathcal{A}_{[n,n+1]}
$$
We aim is to reconstruct a backward quantum Markov chain (see definition \ref{Backw-MCC}),
starting from a sequence $(E_{[n+1,n]})_{n\ge 0})$ of backward Markov transitions expectations.\\
From (\ref{basis2}) the map
\begin{equation}\label{extension_CP}
e_{i_0}^{(0)}e_{i_1}^{(1)}\cdots e_{i_{n-1}}^{(n-1)}e_{i_{n}}^{(n)}e_{i_{n+1}}^{(n+1)} \mapsto e_{i_0}^{(0)}e_{i_1}^{(1)}\cdots e_{i_{n-1}}^{(n-1)}E_{[n+1,n]}(e_{i_{n}}^{(n)}e_{i_{n+1}}^{(n+1)})
\end{equation}
$(i_0,\cdots,i_{n+1})\in I_0\times \cdots \times  I_{n+1}$\\
extends $E_{[n+1,n]}$ to a unique linear map $\widetilde{E}_{n]}$
from $\mathcal A_{n+1]}$ into $\mathcal A_{n]}$.
\begin{lemma}\label{*map}{\rm
$\widetilde{E}_{n]}$ is a $*$-map if and only if it satisfies the following trace--like property
\begin{equation}\label{star_condition}
\widetilde{E}_{n]}(ab)=\widetilde{E}_{n]}(ba)\;  \quad \quad a \in \mathcal A_{n-1]},\  b \in \mathcal{A}_{[n,n+1]}
\end{equation}
}\end{lemma}
\textbf{Proof}. {\rm  For $a \in \mathcal A_{n-1]}$, $b \in \mathcal{A}_{[n,n+1]}$.
If $\widetilde{E}_{n]}$ is a $*$-map then
$$
\widetilde{E}_{n]}(b^* a^*)= \left(\widetilde{E}_{n]}(ab)\right)^*
=\left(a E_{[n,n+1]}(b)\right)^*
=\left(E_{[n,n+1]}(b)\right)^*a^*
$$
Using the completely positivity of $E_{[n,n+1]}$ and the Markov property (\ref{MP}),\\
one gets
$$
\widetilde E_{n]}(b^* a^*)=\widetilde E_{n]}( a^*b^*)
$$
Therefore, $\widetilde{E}_{n]}(ba)=\widetilde{E}_{n]}(ab)$, for each $a \in \mathcal A_{n-1]}$, $b \in \mathcal{A}_{[n,n+1]}$.\\
}
\newpage
\begin{lemma}\label{*map-condition}{\rm The following assertions hold true.
\begin{enumerate}
\item If $\widetilde{E}_{n]}\circ \widetilde E_{n]}^{0}=\widetilde{E}_{n]}$ then $\widetilde{E}_{n]}$ is a $*$-map.
\item $\widetilde{E}_{n]}\circ \widetilde{E}_{n]}^{0}=\widetilde{E}_{n]}$ if and only if  $E_{[n+1,n]}\circ E_{[n+1,n]}^{0}=E_{[n+1,n]}$.
\end{enumerate}
}\end{lemma}
\textbf{Proof}. {\rm \begin{enumerate}
\item For $a \in \mathcal A_{n-1]}$, $b \in \mathcal{A}_{[n,n+1]}$.
If $\widetilde{E}_{n]}\circ \widetilde E_{n]}^{0}=\widetilde{E}_{n]}$ then
$$
\widetilde{E}_{n]}(ab)=\widetilde{E}_{n]}\left(\widetilde E_{n]}^{0}(ab)\right)
=\widetilde{E}_{n]}\left(\widetilde E_{n]}^{0}(ba)\right)
=\widetilde{E}_{n]}(ba)
$$
\item
$
\widetilde{E}_{n]}(ab)=\widetilde{E}_{n]}\left(\widetilde E_{n]}^{0}(ab)\right)
=\widetilde{E}_{n]}\left(a E_{[n+1,n]}^{0}(b)\right)
=aE_{[n+1,n]}\circ E_{[n+1,n]}^{0}(b).
$
\end{enumerate}
From now on we assume that
\begin{equation}\label{star-condition1}
E_{[n+1,n]}\circ E_{[n+1,n]}^{0}=E_{[n+1,n]}
\end{equation}
Therefore $\widetilde{E}_{n]}$ is a $*$--map.\\
\textbf{Remark}. {\rm From (\ref{extension_CP}) the range of $\widetilde{E}_{n]}$ satisfies
$$
\hbox{Range}(\widetilde{E}_{n]} ) \subseteq  \mathcal A_{n-1]}\bigvee \hbox{Range}(E_{[n+1,n]})
$$
and using the Markovianity of $E_{[n+1,n]}$  (see (\ref{MP}) ) one gets:
\begin{equation}\label{range1}
\hbox{Range}(\widetilde{E}_{n]})  \subseteq  \mathcal A_{n-1]}\bigvee \left(\mathcal A_{n-1]}'\cap \mathcal A_n\right)
\end{equation}
}
\begin{theorem}\label{general-QCE}{\rm
The map $\widetilde E_{n]}$ defined through (\ref{extension_CP}) is a Markov quasi-conditional expectation with respect to the following triplet
\begin{equation}\label{reconstruction_triple}
\mathcal {A}_{n-1]}\subseteq \mathcal{A}_{n-1]} \vee  \bigl(\mathcal{A}_{n-1]}^{'}\cap \mathcal {A}_{n}\bigr) \subseteq \mathcal A_{n+1]}
\end{equation}
}
\end{theorem}
\textbf{Proof}. {\rm  By construction and the equation (\ref{star-condition1}) the map $\widetilde{E}_{n]}$ is linear-$*$--map.
Let now move to its complete positivity.\\
For $m \in \mathbb{N}$, let $a_{n],1},\cdots, a_{n],m} \in \mathcal{A}_{n-1]} \vee \bigl({A}_{n-1]}^{'} \cap \mathcal {A}_{n}\bigr)$
and $ a_{n+1],1},\cdots, a_{n+1],m} \in \mathcal A_{n+1]}$.
From (\ref{basis2})  it is enough to consider product elements of the following form
$$
a_{n],i}=a_{n-1],i} a_{n,i},\ \  a_{n-1],i} \in \mathcal{A}_{n-1]}, \  \ a_{n,i}\in ({A}_{n-1]}^{'} \cap \mathcal {A}_{n} ) ,\ \ \ i=1,\cdots, m
$$
$$
a_{n+1],i}=b_{n-1],i} b_{[n,n+1],i}, \ \ b_{n-1],i} \in \mathcal{A}_{n-1]},\ \ b_{[n,n+1],i} \in \mathcal{A}_{[n,n+1]}, \ \ i=1,\cdots, m
$$
\begin{equation}\label{CP}
\sum^m_{j,k=1} a_{n],j}\tilde E_{n]}(a_{n+1],j}  a_{n+1],k}^{*})a_{n],k}^{*}
\end{equation}
$$
=\sum^m_{j,k=1}a_{n-1],j} a_{n,j}\widetilde E_{n]} (b_{n-1],j} b_{[n,n+1],j}b_{[n,n+1],k}^{*} b_{n-1],k}^{*})a_{n,k}^{*} a_{n-1],k}^{*}
$$
One has
$$
\widetilde E_{n]} (b_{n-1],j} b_{[n,n+1],j}b_{[n,n+1],k}^{*} b_{n-1],k}^{*})=b_{n-1],j}E( b_{[n,n+1],j}b_{[n,n+1],k}^{*})b_{n-1],k}^{*}
$$
Then (\ref{CP}) becomes
$$
\sum^m_{j,k=1} a_{n],j}\widetilde E_{n]}(a_{n+1],j}  a_{n+1],k}^{*})a_{n],k}^{*}
=\sum^m_{j,k=1}a_{n-1],j} a_{n,j} b_{n-1],j}E_{[n+1,n]}( b_{[n,n+1],j}b_{[n,n+1],k}^{*})b_{n-1],k}^{*}a_{n,k}^{*} a_{n-1],k}^{*}
$$
$$
=\sum^m_{j,k=1} a_{n,j}E_{[n+1,n]}( b_{[n,n+1],j}b_{[n,n+1],k}^{*})a_{n,k}^{*}(a_{n-1],j}  b_{n-1],j})(a_{n-1],k}b_{n-1],k} )^{*}
$$
Now consider
$$
A= [a_{n,j}E_{[n+1,n]}( b_{[n,n+1],j}b_{[n,n+1],k}^{*})a_{n,k}^{*}] \in  M_m(\mathcal A_{n})
$$
and
$$
B= [(a_{n-1],j}  b_{n-1],j})(a_{n-1],k}b_{n-1],k} )^{*} ] \in  M_m(\mathcal A_{n-1]})
$$
One can check that the matrices $A$ and $B$ are positive.
Then by lemma \ref{schur}, the matrix $C=A\circ B\in  M_m(\mathcal A_{n]})$ is positive.
Therefore, $\widetilde {E}_{n]}$ is completely positive. This complete the prove.\\
}\\
\textbf{Reconstruction of the boundary conditions}. {\rm
For each $n \in \mathbb N $, define
$$
I_{n}=E_{[1,0]}(E_{[2,1]}(\cdots E_{[n+1,n]}(\mathcal A_{n+1}^{+}) \subseteq \mathcal{A}_{0}^{+}
$$
One Remarks that, if all the transition expectations $E_{[n+1,n]}$ are normalized  then $1\in \bigcap_{n\in\mathbb N}I_n$.
}
\begin{lemma}\label{boundary-cd}{\rm
If $J_0 =:\bigcap_{n\geq 0}I_n\ne \{0\}$, then there exist a sequence of $(b_n)_{n\ge 0}$ boundary conditions with respect to the quasi--conditional expectation $(\widetilde {E}_{n]})_{n \geq 0}$.
}\end{lemma}
\textbf{Proof}. {\rm Let fix $b_0\in J_0\setminus \{0\}$ and
define
$$
J_{k}=\bigcap_{n\geq k}I_{n}\quad ;\quad k\in \mathbb N
$$
One can see that
\begin{equation}\label{MC}
\widetilde {E}_{n]}(J_{n+1})=E_{[n+1,n]}(J_{n+1}) =J_n
\end{equation}
Therefore, from \eqref{MC}, we can define a sequence $(b_n)_{n\geq0}\subset  J_n\setminus \{0\}$ satisfying for each $n \in\mathbb N $
\begin{equation}\label{Mtgle-eq}
\widetilde {E}_{n]}(b_{n+1})=b_{n}
\end{equation}
In addition, from the Markov property (\ref{MP}), one has
\begin{equation}\label{MP1}
b_{n}\in \mathcal A_{n}^{+}\cap \mathcal A_{n-1]}^{'}
\end{equation}
Then (\ref{Mtgle-eq}) and $\ref{MP1}$ implies that $(b_{n})_{n\geq}$ is a
sequence of boundary conditions with respect to the sequence  $(\widetilde {E}_{n]})_{n \geq 0}$.\\
}
\textbf{Initial state}. {\rm Let $\phi\in S(\mathcal {A}_{0})$ such that $\phi(b_0)\ne 0$ and define
\begin{equation}\label{initial_state}
\varphi_{0}(a): =\frac{1}{\phi(b_{0})}\phi(a), \ \ \hbox{ for each } a \in \mathcal {A}_{0}
\end{equation}
}
\begin{theorem}\label{recons-BMC}
Under the same conditions as theorem (\ref{general-QCE}) and lemma (\ref{boundary-cd}), the triplet $\{\varphi_{0}, (\widetilde {E}_{n]})_{n \geq 0},(b_{n})_{n\geq0}\}$  defines a backward Markov chain $\varphi$ on $\mathcal A$.
\end{theorem}
\textbf{Proof}. {\rm By construction the triplet $\{\varphi_{0}, (\widetilde {E}_{n]})_{n \geq 0},(b_{n})_{n\geq0}\}$  given respectively by  (\ref{initial_state}), (\ref{Mtgle-eq}) and (\ref{extension_CP}) satisfies the sufficient conditions of Theorem \ref{Suff_Cond-thm}. Then the result follows immediately.}

\section {Examples}

\subsection{Tensor case}

Let $\mathcal M$ be $q\times q$ matrix algebra on $\mathbb{C}$, denote $\mathcal{A}=\bigotimes_{\mathbb{N}}\mathcal{M}$
the tensor product of $\mathbb N$ copies of $\mathcal M$, $j_k : \mathcal M \mapsto j_k(\mathcal M)\subset \mathcal A$
the natural immersion of $\mathcal M$ onto the "$k-th$ factor" of the product $\bigotimes_{\mathbb{N}}\mathcal{M}$ and
$\mathcal A_{[m;n]}$ the C$^*$sub-algebra of $\mathcal A$ spanned by $\bigcup_{k=m}^{n}j_k(\mathcal M).$

\begin{theorem}
Let $E_{[n+1,n]} :\mathcal{A}_{[n,n+1]} \to \mathcal{A}_{\{n\}}$ be a completely positive linear map.
The formula
\begin{equation}\label{extension_CP1}
a\otimes b \mapsto a\otimes E_{[n+1,n]}(b) \quad ; \quad a\in\mathcal{A}_{n-1]},\  \, b\in \mathcal{A}_{[n,n+1]}
\end{equation}
determines a unique  quasi--conditional expectation $\widetilde{E}_{n]}$ with respect to the triplet $\mathcal{A}_{n-1]}\subseteq
\mathcal{A}_{n]} \subseteq \mathcal{A}_{n+1]}$.
\end{theorem}
\textbf{Proof}.
By linearity it is enough to prove complete positivity for elements of the form
$$
x_i = u_i\otimes v_i\in\mathcal{A}_{n-1]} \otimes \mathcal {A}_{n},\quad  y_i
=s_i\otimes t_i \in \mathcal{A}_{n-1]}\otimes \mathcal{A}_{[n,n+1]}
$$
One has
\begin{equation}\label{reforme_cp}
\sum_{i,k}x_i^\ast\widetilde{E}_{n]}(y_i^\ast y_k)x_k
=\sum_{i,k}(u_i^\ast\otimes v_i^\ast)\widetilde{E}_{n]}\left(s_i^\ast\otimes t_i^\ast s_k\otimes t_k\right)
(u_k\otimes v_k) \\
\end{equation}
$$
=\sum_{i,k}(u_i^\ast\otimes v_i^\ast)
\widetilde{E}_{n]}\left(s_i^\ast s_k\otimes t_i^\ast  t_k\right)(u_k\otimes v_k) \\
=\sum_{i,k}u_i^\ast s_i^\ast s_ku_k\otimes v_i^\ast E_{[n+1,n]}(t_i^\ast  t_k)v_k
$$
Now consider $A= (u_i^\ast s_i^\ast s_k u_k)\in  M_n(\mathcal{A}_{n-1]})$ and
$B= (v_i^\ast E_{[n+1,n]}(t_i^\ast  t_k) v_k) \in  M_n(\mathcal{A}_{n})$
For $a= (a_1, \cdots, a_n)^T\in  \mathcal{A}_{n-1]}^n$
$$
a^\ast A a = \sum_{i,k}a_i^\ast u_i^\ast s_i^\ast s_k u_ka_k = |\sum_{i}s_i u_ia_i|^2
$$
therefore $A\in  M_n(\mathcal{A}_{n-1]})^{+}$.
Similarly, for $b= (b_1, \cdots, b_n)^T\in  \mathcal{A}_{n}^n$, taking in account the complete
positivity of $E_{[n+1,n]}$ one gets
$$
b^\ast B b = \sum_{i,k}b_i^\ast v_i^\ast E_{[n+1,n]}(t_i^\ast  t_k) v_kb_k
=  \sum_{i,k} (v_ib_i)^\ast E_{[n+1,n]}(t_i^\ast  t_k) (v_kb_k)
\geq 0
$$
therefore  $B\in  M_n(\mathcal{A}_{n})^{+}$.
From lemma \ref{Tensor_Schur_lemma} one then gets
$A\circ^\otimes B\in  M_n(\mathcal{A}_{n]})^+$. In particular, denoting
${\bf 1}_{n,\mathcal{A}_{n]}}
:= \left( 1_{\mathcal{A}_{n]}}, \cdots, 1_{\mathcal{A}_{n]}}\right)^T\in (\mathcal{A}_{n]})^n$, one has
$$
\sum_{i,k}\left(u_i^\ast s_i^\ast s_k u_k\right)\otimes \left( v_i^\ast E(t_i^\ast  t_k) v_k\right) = {\bf 1}_{n, \mathcal{A}_{n]}}^T A\circ^\otimes B{\bf 1}_{n, \mathcal{A}_{n]}} \ge 0
$$
i.e. the right hand side of equation  (\ref{reforme_cp}) is positive and this ends the proof.\\

\textbf{Remark}.\rm{ In the this case $\mathcal{A}_{n-1]}^{'}\bigcap \mathcal{A}_{[n,n+1]}=\mathcal{A}_{[n,n+1]}$, which means that, the Umegaki conditionals expectations  $E_{[n+1,n]}^{0}$ are the identity of $\mathcal{A}_{[n,n+1]}$.
}
\subsection{Fermi case}\label{Fermi-case}

In this section  $\mathcal A$  is the Fermi algebra generated by a family of creators and annihilators
$\{ a_{i}, a_{i}^+\quad; i\in \mathbb N\}$ and relations
\begin{equation}\label{anti_commutation}
(a_j)^+ =a_j^+\ , \quad \{a^+_j, a_k\} = \delta_{jk }1_{\mathcal A}\ , \quad \{a_j, a_k\} =0, \ j,k \in I
\end{equation}
For any $J\subseteq \mathbb N$, denote $\mathcal A(J)$ the sub-algebra generated by  $\{ a_j, a_j^+,\quad j\in J\}$.
Now consider any partition  $(J_n)_{n\in \mathbb N} $ of the set $\mathbb N$ such that for each $n$ the set $J_n$ is finite.
Put $d_n = |J_n|<\infty$.\\
Let $\mathcal A_{n } = \mathcal A(J_n)$, it is then the Fermi subalgebra of $\mathcal A$
generated by the $2d_n$ generators $a_{1}, a_{1}^+,\, \cdots\, a_{d_n}, \, a_{d_n}^+$.\\
In this notations one gets for each $I\subseteq \mathbb N$,
$$
\mathcal A_I = \bigvee_{n\in I} \mathcal A_{n} = \bigvee_{n\in I} \mathcal A(J_n) = \mathcal A(\bigcup_{n\in I}J_n)
$$
In particular
$$
\mathcal A_{n]} = \mathcal A(J_0\cup\cdots J_n)
$$
Let $J \subset \mathbb N$ finite and let $m=|J|$. For each $j\in J$ the elements $a_j, a_{j}^{+}, a_ja_{j}^+, a_{j}^+a_{j}$
for a linear basis of the sub-algebra $\mathcal A(\{j\})$ generated by $a_j$ and $ a_j^+$.\\
Since $a_{j}^{*}$ and $a_{j}$ anti--commute among different indices, $a_{j}^{*}$ and $a_{j}$
with a specific $j$ can be brought together at any spot in a monomial, with possible sign change
(without changing the ordering among themselves), and this can be done for each $j$.
Fix an enumeration $i_1, i_2,\dots, i_m$ of $J$. Therefore, the monomials of the form
\begin{equation}\label{monome_fermi}
A_{i_1}A_{i_2}\cdots A_{i_{m}}
\end{equation}
where $A_{j}$ is one of $a_j, a_{j}^{+}, a_ja_{j}^+, a_{j}^+a_{j}$,
consists a linearly spanning family of cardinality $4^{m}$.\\
In the other hand, the  Jordan-Klein-Wigner transformation establishes the (linear) isomorphism
\begin{equation}\label{JKV_isomorphism}
\mathcal A(J )\sim \bigotimes_{J}M_2(\mathbb C)
\end{equation}
In fact, denote
$$ 1= \left(\begin{array}{cc}
                           1 & 0 \\
                           0 & 1
                         \end{array}\right),\quad\sigma_z = \left(\begin{array}{cc}
                           -1 & 0 \\
                           0 & 1
                         \end{array}\right) $$

Put for each $j\in [1, m]$
\begin{equation}\label{JKV_basis}
e_{kl}(j) =  \underbrace{\sigma_z\otimes \cdots \sigma_z }_{j-1\;  times}\otimes e_{kl}\otimes \underbrace{1 \otimes \cdots \otimes 1}_{m-j\ times}, \quad 1\le k,l\le 2
\end{equation}
 where $(e_{kl})_{1\le k,l\le 2}$ is the canonical system of  $M_2(\mathbb C)$.
Then the identification
$$
a_{i_j} \mapsto e_{21}(j)\quad ; \quad a_{i_j}^+ \mapsto e_{12}(j)
$$
$$
a_{i_j}^+a_{i_j} \mapsto e_{11}(j)\quad ; \quad a_{i_j}a_{i_j}^+ \mapsto e_{22}(j)
$$
  realizes the isomorphism (\ref{JKV_isomorphism}).
Therefore, monomials (\ref{monome_fermi}) consist a linear basis of the sub-algebra $\mathcal A(J)$.\\

\begin{definition}{\rm
$\Theta_{J}$ denotes the unique automorphism of $\mathcal{A}$ satisfying
\begin{equation}\label{Theta}
\Theta_{J}(a_{i})=-a_{i},  \quad \Theta_{J}(a_{i}^{+})=-a_{i}^{+}, \quad (i \in J)
\end{equation}
$$
\Theta_{J}(a_{i})=a_{i},  \quad \Theta_{J}(a_{i}^{*})=a_{i}^{+}, \quad (i \in J^{c})
$$
In particular, we denote $\Theta= \Theta_{\mathbb N}$.
}\end{definition}
The even and odd parts of $\mathcal{A}$ are defined as
\begin{equation}
\mathcal{A}_{+}\equiv\{a \in \mathcal{A}\mid \Theta(a)=a \}, \ \  \mathcal{A}_{-}\equiv\{a \in \mathcal{A}\mid \Theta(a)=-a \}.
\end{equation}
\textbf{Remark}. {\rm
Such $\Theta$ exists and is unique because (\ref{Theta}) preserves CAR. It obviously satisfies
$$
\Theta^{2}=1
$$
}
\textbf{Remark}. {\rm
For any $a \in \mathcal{A}_{J}$
$$
a=a_{+}+ a_{-},\ \ a_{\pm}=\frac{1}{2}(a \pm \Theta(a))
$$
gives the (unique) splitting of $a$ into a sum of $a_{+}\in \mathcal{A}_{\{J,+\}}$ and $a_{-}\in \mathcal{A}_{\{J,-\}}$,
where the even and odd parts of $\mathcal{A}_{J}$ are denoted by $\mathcal{A}_{\{J,+\}}$ and $\mathcal{A}_{\{J,-\}}$.
}
\begin{definition}\label{df-even}{\rm
A map $E$: $\mathcal A$ $\rightarrow$ $\mathcal B$ between the Fermi algebra $\mathcal A$, $\mathcal B$
is said to be even if
$$
E\circ \Theta= E
$$
}\end{definition}
\textbf{Remark}. {\rm
If $E$ is even then for each $a \in \mathcal {A}_{-}$
\begin{equation}\label{even-prt}
E(a)=E(\Theta(a))=-E(a)=0,
\end{equation}
}
\begin{lemma}\label{commutant}{\rm
For a finite $J \in \mathbb N$,
\begin{equation}
(\mathcal{A}_{J})^{'} = \mathcal{A}_{\{J^{c},+\}} + v_{J} \mathcal{A}_{\{J^{c},-\}},
\end{equation}
where $v_J $ is the self-adjoint unitary in $\mathcal{A}_{\{J,+\}}$  given by
\begin{equation}
v_{J}\equiv\prod_{n\in J}v_{n},\ \ v_{n}=\prod_{i=1}^{d_n}a_{i}a_{i}^{+}-a_{i}^{+}a_{i}
\end{equation}
}\end{lemma}
\textbf{Proof}.(see \cite{[ArMo]})\\
\textbf{Remark}. {\rm By lemma \ref{commutant}, the Umegaki conditional expectation $E^{0}_{[n+1,n]}$ are defined by
$$
E^{0}_{[n+1,n]}:\ \mathcal{A}_{[n,n+1]} \mapsto \mathcal{A}_{\{[n,n+1],+\}}
$$
$$
b=b_{+}+ b_{-} \mapsto b_{+}
$$
}
\begin{lemma}\label{even-cond}{\rm
$E_{[n+1,n]}$ is even if and only if  $E_{[n+1,n]}\circ E_{[n+1,n]}^{0}=E_{[n+1,n]}$.
}\end{lemma}
\textbf{Proof}. {\rm Let $b \in \mathcal{A}_{[n,n+1]}$.\\
If $E_{[n+1,n]}$ is even then the unique splitting of $b$ into a sum of $b_{\pm} \in \mathcal{A}_{[n,n+1],\pm}$
implies that
$$
E_{[n+1,n]}(b)=E_{[n+1,n]}(\Theta(b))=E_{[n+1,n]}(b_{+}-b_{-})=E_{[n+1,n]}(b_{+})=E_{[n+1,n]}\left( E_{[n+1,n]}^{0}(b)\right)
$$
}
\begin{theorem}\label{fermi-QCE}{\rm Let
$E_{[n+1,n]}$ be a even backward Markov transition expectation from
$\mathcal {A}_{[n,n+1]}\rightarrow \mathcal {A}_{n}$, then the map $\widetilde E_{n]}$
defined through (\ref{extension_CP}) is a quasi--conditional expectation with respect to the triplet\\
$\mathcal {A}_{n-1]}\subseteq \mathcal {A}_{n-1]} \vee \mathcal {A}_{\{n,+\}} \subseteq \mathcal {A}_{n+1]}$.
}\end{theorem}
\textbf{Proof}. {\rm
From (\ref{extension_CP}), $\widetilde E_{n]}$ is a linear map.\\
For $a \in \mathcal{A}_{n-1]}$ and $b\in \mathcal{A}_{[n,n+1]}$, we have
$$
\widetilde E_{n]}(ab)=aE_{[n+1,n]}(b)=aE_{[n+1,n]}(b_{+}+b_{-})=a E_{[n+1,n]}(b_{+})
$$
And
$$
\widetilde E_{n]}(ba)=\widetilde E_{n]}((b_{+}+b_{-})(a_{+}+a_{-}))=\widetilde E_{n]}(b_{+}a_{+}+b_{-}a_{+}+b_{+}a_{-}+b_{-}a_{-})
$$
Since $a \in \mathcal{A}_{n-1]}$ and $b\in \mathcal{A}_{[n,n+1]}$, we have
$$
\widetilde E_{n]}(ba)=\widetilde E_{n]}(a_{+}b_{+}+a_{+}b_{-}+a_{-}b_{+}-a_{-}b_{-})
$$
By linearity of $\widetilde E_{n]}$, we get
$$
\widetilde E_{n]}(ba)=\widetilde E_{n]}(a_{+}b_{+})+\widetilde E_{n]}(a_{+}b_{-})+\widetilde E_{n]}(a_{-}b_{+})-\widetilde E_{n]}(a_{-}b_{-})
$$
$$
=a_{+}E_{[n+1,n]}(b_{+})+a_{+}E_{[n+1,n]}(b_{-})+a_{-}E_{[n+1,n]}(b_{+})-a_{-}E_{[n+1,n]}(b_{-})
$$
Since $E_{[n+1,n]}$ is even, we obtain
$$
\widetilde E_{n]}(ba)=a_{+}E_{[n+1,n]}(b_{+})+a_{-}E_{[n+1,n]}(b_{+})=aE_{[n+1,n]}(b_{+})=\widetilde E_{n]}(ab)
$$
Therefore, $\widetilde E_{n]}$  satisfy the trace--like property. Then by lemma (\ref{*map}),  $\widetilde E_{n]}$ define a linear $*$--map.
And yet, from lemma (\ref{commutant}) and (\ref{range1}), one has
$$
\hbox{Range}(\widetilde{E}_{n]})  \subseteq  \mathcal A_{n-1]}\bigvee \mathcal {A}_{\{n,+\}}
$$
Let now move to its complete positivity.\\
For $m \in \mathbb{N}$, let $a_{n],1},\cdots,a_{n],m}\in \mathcal{A}_{n-1]} \vee \mathcal{A}_{\{n,+\}}$ and $a_{n+1],1},\cdots a_{n+1],m} \in \mathcal A_{n+1]}$.
From (\ref{basis2}) we can rewrite those elements in the following form
$$
a_{n],i}=a_{n-1],i} a_{n,i},\ \  a_{n-1],i} \in \mathcal{A}_{n-1]}, \  \ a_{n,i}\in \mathcal{A}_{\{n,+\}},\ \ \ i=1,\cdots, m
$$
$$
a_{n+1],i}=b_{n-1],i} b_{[n,n+1],i}, \ \ b_{n-1],i} \in \mathcal{A}_{n-1]},\ \ b_{[n,n+1],i} \in \mathcal{A}_{[n,n+1]}, \ \ i=1,\cdots, m
$$
\begin{equation}\label{CP}
\sum^m_{j,k=1} a_{n],j}\widetilde E_{n]}(a_{n+1],j}  a_{n+1],k}^{+})a_{n],k}^{+}
\end{equation}
$$
=\sum^m_{j,k=1}a_{n-1],j} a_{n,j}\widetilde E_{n]} (b_{n-1],j} b_{[n,n+1],j}b_{[n,n+1],k}^{+} b_{n-1],k}^{+})a_{n,k}^{+} a_{n-1],k}^{+}
$$
One has
$$
\widetilde E_{n]} (b_{n-1],j} b_{[n,n+1],j}b_{[n,n+1],k}^{+} b_{n-1],k}^{+})=b_{n-1],j}E_{[n+1,n]}( b_{[n,n+1],j}b_{[n,n+1],k}^{+})b_{n-1],k}^{+}
$$
Then (\ref{CP}) becomes
$$
\sum^m_{j,k=1} a_{n],j}\widetilde E_{n]}(a_{n+1],j}  a_{n+1],k}^{+})a_{n],k}^{+}
$$
$$
=\sum^m_{j,k=1}a_{n-1],j} a_{n,j} b_{n-1],j}E_{[n+1,n]}( b_{[n,n+1],j}b_{[n,n+1],k}^{+})b_{n-1],k}^{+}a_{n,k}^{+} a_{n-1],k}^{+}
$$
$$
=\sum^m_{j,k=1} a_{n,j}E_{[n+1,n]}( b_{[n,n+1],j}b_{[n,n+1],k}^{+})a_{n,k}^{+}(a_{n-1],j}  b_{n-1],j})(a_{n-1],k}b_{n-1],k} )^{+}
$$
Now consider
$$
A= [a_{n,j}E_{[n+1,n]}( b_{[n,n+1],j}b_{[n,n+1],k}^{+})a_{n,k}^{+}]\in  M_m(A_{n]})
$$
and
$$
B= [(a_{n-1],j}  b_{n-1],j})(a_{n-1],k}b_{n-1],k} )^{+} ] \in  M_m(\mathcal A_{n]})
$$
One can see that the matrices $A$ and $B$ are positive. Then by lemma \ref{schur}, the matrix
$$
C=A\circ B=[a_{n,j}E_{[n+1,n]}( b_{[n,n+1],j}b_{[n,n+1],k}^{+})a_{n,k}^{+}(a_{n-1],j}  b_{n-1],j})(a_{n-1],k}b_{n-1],k} )^{+}]
$$
is positive. Therefore, we obtain that $\widetilde E_{n]}$ is completely positive.}
\subsection*{Acknowledgment}

\end{document}